\documentclass[11pt]{amsart}

\usepackage{amsmath, amsthm, amssymb, amsfonts, mathrsfs, amscd, amsthm}
\usepackage[all]{xy}

\newtheorem{thm}{Theorem}[section]
\newtheorem{pro}[thm]{Proposition}
\newtheorem{lem}[thm]{Lemma}
\newtheorem{cor}[thm]{Corollary}

\theoremstyle{definition}
\newtheorem{de}[thm]{Definition}

\theoremstyle{remark}
\newcommand{\rem}{{\it Remark}}

\newcommand{\mc}[1]{\mathcal{#1}}

\newcommand{\mf}[1]{\mathfrak{#1}}

\newcommand{\tr}[1]{\textrm{#1}}

\newcommand{\defi}{\begin{de}}
\newcommand{\ed}{\end{de}}
\newcommand{\bl}{\begin{lem}}
\newcommand{\el}{\end{lem}}
\newcommand{\bp}{\begin{pro}}
\newcommand{\ep}{\end{pro}}
\newcommand{\bt}{\begin{thm}}
\newcommand{\et}{\end{thm}}
\newcommand{\bc}{\begin{cor}}
\newcommand{\ec}{\end{cor}}
\newcommand{\bpf}{\begin{proof}}
\newcommand{\epf}{\end{proof}}

\newcommand{\beq}{\begin{equation}}
\newcommand{\eeq}{\end{equation}}
\newcommand{\beqs}{\begin{equation*}}
\newcommand{\eeqs}{\end{equation*}}
\newcommand{\ben}{\begin{enumerate}}
\newcommand{\een}{\end{enumerate}}
\newcommand{\bit}{\begin{itemize}}
\newcommand{\eit}{\end{itemize}}

\newcommand{\ul}[1]{\underline{#1}}

\begin {document}

\title{Projective Varieties covered by isotrivial families }

\author{Anupam Bhatnagar} 



\address{Department of Mathematics, Borough of Manhattan Community College, The City University of New York; 199 Chambers Street, New York, NY 10007, U.S.A.}

\email{anupambhatnagar@gmail.com}

\begin{abstract}
Let $X,Y$ be projective schemes over a discrete valuation ring $R$, where $Y$ is generically smooth and $g: X \to Y$ a surjective $R$-morphism such that  $g_*\mc{O}_{X} = \mc{O}_{Y}$. We show that if the family $X \to Spec(R)$ is isotrivial, then the generic fiber of the family $Y\to Spec(R)$ is isotrivial.
\end{abstract}

\maketitle

\section{Introduction}

Let $k$ be a field of characteristic zero, $C$ a smooth projective curve defined over $k$ with function field $F$. 

\defi Let $S$ be a scheme over $k$ and $\pi: X\to S$ a flat family of schemes. The family $\pi$ is called {\it trivial} if there exists a scheme $X_0$ defined over $k$ such that $X \cong X_0 \times_k S$ and it is called {\it isotrivial}
if there exists a finite surjective \'etale extension $S' \to S$ such that $\pi_{S'}: X \times_S S' \to S'$ is trivial. 
\ed

Let $R$ be the local ring at a closed point $P\in C$. Let $X,Y$ be projective schemes over $Spec(R)$ where $Y$ is generically smooth and $g: X\to Y$ a surjective $R$-morphism such that $g_*\mc{O}_X = \mc{O}_Y$. We show in Theorem \ref{isotrivial thm}, if the family $X\to Spec(R)$ is isotrivial, then the generic fiber of the family $Y \to Spec(R)$ is isotrivial. 

\smallskip

\noindent{\it{Sketch of proof}}. The deformations of $Y$ are controlled by its differentials. We 
study the deformations locally i.e. over the discrete valuation ring $R$, and show that the
fundamental exact sequence associated to the diagram  
$$X \to Y \stackrel{p}{\to} Spec (R) \to Spec(k)$$
\beq 
0 \to p^* \Omega_{Spec(R)/Spec(k)} \to \Omega_{Y/Spec(k)} \to \Omega_{Y/Spec(R)} \to 0
\eeq
is split exact and consequently the deformations of $Y$ are governed by that of $X$. 
Then we consider an infinitesimal deformation of $Y \to Spec(R)$ over the henselization of $R$ (denoted by 
$\tilde{R}$) and show that the sequence above remains split exact at every level of the deformation. Finally we 
use a result of Greenberg to pass from $\tilde{R}$ to $R$. 

The proof can probably be adapted as is to the positive characteristic case if we assume the extension $F/k$ is separable. As an application we shall use this result in a forthcoming paper with Lucien Szpiro, on parametrization of points of canonical height zero of an algebraic dynamical system. 

In \cite{CH1}, (\cite{CH2}, Thm. 3.3, p. 702) Chatzidakis and Hrushovski answer the same question using model 
theoretic methods. Their methods are intrinsically birational thus they have a slightly less precise conclusion. They assume the extension $F/k$ is regular. If the ground field $k$ is perfect and $Y$ is one dimensional then their result extends to positive characteristic.

{\it Notation}. 
Throughout this paper $k$ denotes a field of characteristic zero, $C$ a smooth projective curve over $k$ with function field $F$. Given a scheme $X$ over $C$ we denote its generic fiber by $X_F$.


\section{Descent}

In this section we assume that $k$ is algebraically closed. The proofs work without this hypothesis with some minor modifications. Let $P$ be a closed point on the curve $C$ and $R=\mc{O}_{C,P}$, the local ring at $P$.
For ease of notation we denote Spec($R$) and Spec($k$) by $R$ and $k$ respectively in the sheaves of differentials. 

\bp \label{diff thm}
Let 
$$ X \stackrel{g}{\to} Y \stackrel{p}{\to} Spec(R) \to Spec(k)$$ 
be morphisms of schemes where $X,Y$ are projective, $X\to Spec(R)$ is an isotrivial family, $Y$ is reduced and 
$g:X \to Y$ a surjective $R$-morphism such that $g_*\mc{O}_{X}= \mc{O}_{Y}$. Then the sequence of differentials 
on $Y$,
\beq \label{diff seq}
p^* \Omega_{R/k} \to \Omega_{Y/k} \to \Omega_{Y/R} \to 0 \eeq
is split exact
\ep
\bpf After a quasi-finite unramified base change $Spec(R') \to Spec(R)$ we may assume that the family 
$X_{R'} \to Spec(R')$ is trivial. Slightly abusing the notation we shall say the family $X \to Spec(R)$ is trivial. It follows that the sequence of differentials on $X$ 
\beq \label{seq1}
0 \to g^*p^*\Omega_{R/k} \to \Omega_{X/k} \to \Omega_{X/R} \to 0 \eeq 
is split exact. Since pullbacks preserve right exactness, the sequence of differentials on $Y$ pulled back to $X$ along $g$
\beq \label{seq2}
g^*p^*\Omega_{R/k} \to g^*\Omega_{Y/k} \to g^*\Omega_{Y/R} \to 0 \eeq 
is exact. The morphism $g: X \to Y$ induces the following commutative diagram:
\beq \label{diagram1}
\xymatrix{ & g^*p^*\Omega_{R/k} \ar[r] \ar[d]^{Id} & g^*\Omega_{Y/k} \ar[r] \ar[d] & g^*\Omega_{Y/R} \ar[r] \ar[d] & 0 \\
0 \ar[r] & g^*p^*\Omega_{R/k} \ar[r] & \Omega_{X/k} \ar[r] & \Omega_{X/R} \ar[r] & 0} 
\eeq
It follows that (\ref{seq2}) is split exact. Since $g_*$ preserves direct sums, we have 
$$g_*g^*\Omega_{Y/k} \cong g_*g^*\Omega_{Y/R} \oplus g_*g^*p^*\Omega_{R/k}$$ 
The natural map $\Omega_{Y/k} \to g_*g^*\Omega_{Y/k}$ induces the following commutative diagram:
$$ \xymatrix{ & p^*\Omega_{R/k} \ar[r] \ar[d] & \Omega_{Y/k} \ar[r] \ar[d] & \Omega_{Y/R} \ar[r] \ar[d] & 0 \\
0 \ar[r] & g_*g^*p^*\Omega_{R/k} \ar[r] & g_*g^*\Omega_{Y/k} \ar[r] & g_* g^*\Omega_{Y/R} \ar[r]  & 0
} $$
Note that the bottom row is split exact, $p^*\Omega_{R/k} \cong \mc{O}_Y$ and $g_*g^*p^*\Omega_{R/k} \cong g_*\mc{O}_X$. By assumption $g_*\mc{O}_X = \mc{O}_Y$ thus the top row is split exact. 
\epf

\rem. The condition $g_*\mc{O}_ X = \mc{O}_Y$ implies that $g$ has connected fibers (\cite{H}, Cor. 11.3, p. 279). If $g$ is finite, then the condition $g_*\mc{O}_ X = \mc{O}_Y$ implies that $\deg(g) =1$, hence $g$ is birational. Moreover if $Y$ is normal, using Zariski's Main Theorem (\cite{M}, p. 209) we conclude that $g$ is an isomorphism. 

We now consider an infinitesimal deformation of $Y$ over the henselian discrete valuation 
ring, denoted $\tilde{R}$ and proceed to show that the family $Y_F \to$ Spec($F$) is isotrivial. Before we 
proceed we need the following definitions:

\defi Let $S$ be a smooth scheme of finite type over $k$ and $f: X\to S$ a morphism of schemes. If $f$ is smooth then the sequence
\beq 0 \to f^*\Omega_{S/k} \to \Omega_{X/k} \to \Omega_{X/S} \to 0 
\eeq
is exact. This extension is non-trivial in general and is given by a class $c \in Ext^1(\Omega_{X/S}, f^*\Omega_{S/k})$. 
Since $\Omega_{X/S}$ is locally free, one has 
$$Ext^1_{\mc{O}_X}(\Omega_{X/S}, f^*\Omega_{S/k}) \cong  Ext^1_{\mc{O}_X}(\mc{O}_X, T_{X/S} \otimes f^*\Omega_{S/k})  \cong  H^1(X,T_{X/S} \otimes f^* \Omega_{S/k})  $$
The image of $c$ by the canonical map 
\begin{eqnarray*}
H^1(X,T_{X/S} \otimes f^* \Omega_{S/k}) \to &H^0(S, R^1f_*(T_{X/S} \otimes f^*\Omega_{S/k}))  \\
& || \\
& H^0(S, R^1f_*T_{X/S} \otimes \Omega_{S/k} )
\end{eqnarray*}
is called the {\it Kodaira-Spencer class of $X/S$}. One can view this class as a morphism also i.e. the Kodaira-Spencer morphism 
$$ \kappa_{X/S}: T_S \to R^1f_* T_{X/S} $$ The fiber $(\kappa_{X/S})_s = \kappa_s: T_{S,s} \to H^1(X_s, T_{X_s})$ is the Kodaira-Spencer map at $s\in S$. 
\ed
The Kodaira-Spencer map at $s$ measures how $X_s$ deforms in the family $X/S$ in the neighborhood of $s$ (\cite{B}, p. 165).

\defi $($\cite{D}, p. 255$)$ A local ring $A$ is {\it henselian} if every finite $A$-algebra 
$B$ is a product of local rings. We define the {\it henselization} of $A$ to be a pair 
$(\tilde{A},i)$, where $\tilde{A}$ is a local henselian ring and $i:A\to \tilde{A}$ is 
a local homomorphism such that: 
for any local henselian ring $B$ and any local homomorphism $u:A\to B$ there exists a 
unique local homomorphism $\tilde{u}: \tilde{A} \to B$ such that $u = \tilde{u} \circ i$. 
\ed

From here on we assume that $Y$ is generically smooth. Let $R$ be the local ring at $P\in C$, $\tilde{R}$ its henselization, $\tilde{\mf{m}}$ the maximal ideal of $\tilde{R}$. 
Define $R_n =\tilde{R}/\tilde{\mf{m}}^{n+1} $ for each $n\geq 0$. There are natural maps 
$Spec(\tilde{R}) \to Spec(R)$ and $Spec(R_{n-1}) \to Spec(R_n)$ induced by the projections 
$R_n \to R_{n-1}$ for $n\geq 1$. Define $\tilde{Y}:= Y \times_R Spec(\tilde{R})$, $Y_n:= \tilde{Y} \times_R Spec(R_n)$ for each $n\geq 0$. We have the following commutative diagram of schemes: 
\begin{align}
\xymatrix{ 
Y_F \ar[r] \ar[d] & Y \ar[d]^p & \tilde{Y} \ar[l]  \ar[d] \ldots  & Y_n \ar[l] \ar[d]^{p_n} \ldots &  Y_0 \ar[l] \ar[d]^{p_0} \\
Spec(F) \ar[r] & Spec(R) & Spec(\tilde{R}) \ar[l] \ldots & Spec(R_n) \ar[l] \ldots &Spec(k) \ar[l]} 
\end{align}

\bp \label{KS prop}
For each $n\geq 0$ the sequence of differentials associated to $Y_n \to Spec(R_n)$ i.e.
\beq \label{seq 3}
0 \to p^*_n\Omega_{R_n/k} \to \Omega_{Y_n/k} \to \Omega_{Y_n/R_n} \to 0 \eeq 
is split exact. Moreover, $Y_n \to Spec(R_n)$ is trivial. 
\ep
\bpf Pulling back (\ref{diff seq}) along the natural map $Spec(R_n) \to Spec(R)$ we get the sequence (\ref{seq 3}). Since pullbacks preserve direct sums, the sequence (\ref{seq 3}) is split exact i.e. the Kodaira-Spencer class of  $Y_n /Spec(R_n)$ is trivial. In other words, $Y_n \to Spec(R_n)$ is trivial. It follows that $Y_n \cong Y_0 \times_k Spec(R_n) \tr{ for each } n\geq 0$.
\epf

\defi If $V,W$ and $T$ are $S$-schemes, an {\it $S$-isomorphism from $V$ to $W$} parametrized by $T$ will mean a $T$-isomorphism from $V \times_S T \to W \times_S T$. The set of all such isomorphisms will be denoted by $\ul{Isom}_S(V,W)(T)$. \ed

The association $T \mapsto \ul{Isom}_S(V,W)(T)$ defines a contravariant functor 
$$\ul{Isom}_S(V,W) : (Sch/S)^{\circ} \to (Sets)$$ 
The functor $\ul{Isom}_S(V,W)$ is representable  whenever $V,W$ are flat and projective over $S$. For a proof of the representability of the $\ul{Isom}$ functor we refer the reader to (\cite{FGA} p. 132-133). We denote the scheme representing the functor $\ul{Isom}_S(V,W)$ by $Isom_S(V,W)$. 

To conclude that the family $Y_F \to Spec(F)$ is isotrivial we need the following result of Greenberg:

\bt \label{greenberg thm}Let $\tilde{R}$ be a henselian discrete valuation ring, with $t$ the generator of the 
maximal ideal. Let $\tilde{Z}$ be a scheme of finite type over $\tilde{R}$. Then $\tilde{Z}$ has a point in $\tilde
{R}$ if and only if $\tilde{Z}$ has a point in $\tilde{R}/t^n$ for every $n\geq 1$.
\et
\bpf \cite{G}, Corollary 2. \epf

\bt \label{isotrivial thm}
Let $X,Y$ be projective schemes over a discrete valuation ring $R$ where $X\to Spec(R)$ is an isotrivial family, $Y$ is generically smooth and $g: X\to Y$ a surjective $R$-morphism such that $g_*\mc{O}_X =\mc{O}_Y$. Then $Y\to Spec(R)$ is generically isotrivial i.e. 
$$ Y_{F'} \cong Y_0 \times_k Spec(F')$$ 
where $F'$ is a finite extension of $F$.
\et 
\bpf Observe that $\tilde{Y}$ and $Y_0 \times_k Spec(\tilde{R})$ are flat, projective over $Spec(\tilde{R})$. Let $\ul{Isom}_{\tilde{R}}(\tilde{Y},Y_0 \times_k  Spec(\tilde{R}))(T)$ be the set of isomorphisms from 
$$ \tilde{Y} \times_{\tilde{R}} T \to (Y_0 \times_k Spec (\tilde{R})) \times_{\tilde{R}} T $$
Let $f \in Isom_{\tilde{R}}(\tilde{Y}, Y_0 \times_k Spec (\tilde{R}))(k)$ and $\Gamma_f$, denote the graph of $f$. 
Let $L$ be a very ample line bundle of $Y$ and $\tilde{L}$ the pullback of $L$ to $Y$ along the morphism $\tilde{Y} \to Y$.
Let $P(t)$ be the Hilbert polynomial of $(\Gamma_f)_k$, the special fiber of $\Gamma_f$, with respect to 
$\tilde{L}$. Then the functor
$$ \ul{Isom}_{\tilde{R}}(\tilde{Y}, Y_0 \times_k Spec (\tilde{R}))  \cap \ul{Hilb}^{P(t)}_{\tilde{Y} \times_{\tilde R} (Y_0 \times_k Spec (\tilde{R}))} $$
is representable and the scheme representing it, denoted by $Z'$ is of finite type. For 
$$Z =Isom_R(Y, Y_0 \times Spec(R)) \cap Hilb^{P(t)}_{Y\times_R (Y_0 \times_k Spec(R))}$$ 
we have $Z' = Z \times_R Spec(\tilde{R})$.
By Proposition \ref{KS prop},
$Y_n \cong Y_0 \times_k Spec(R_n)$ for each $n\geq 0$. Thus $Z'$ has a $R_n$-point for every 
$n\geq 0$. By the previous theorem $Z'$ has a $\tilde{R}$-point i.e. 
$\tilde{Y} \cong Y_0 \times_k Spec(\tilde{R})$. Note that $\tilde{R}$ is a limit of etale covers of $R$ so there exists 
an etale cover $R'$ of $R$ such that $Y_{R'} \cong Y_0 \times_k Spec(R')$. Thus $F'$ (the quotient field of $R'$) 
satisfies the requirements of the theorem. 
\epf


{\it Acknowledgements}. I thank Lucien Szpiro for introducing me to the subject of algebraic dynamics and for sharing this question during my PhD. Many thanks to Raymond Hoobler for several interesting conversations towards this paper. Thanks to Antoine Chambert-Loir on giving valuable feedback on an earlier version of this paper and to Madhav Nori for suggestions toward simplifying the proof of Proposition \ref{diff thm}. I also thank the referee who gave me valuable feedback to improve the exposition of this article.


\end{document}